\documentclass{article}

\usepackage[repositionpullbacks,midshaft,PostScript=dvips,nohug,scriptlabels,size=2em]{diagrams}

\usepackage{latexsym}
\usepackage{amssymb}

\newcommand{\dashbk}{-}
\newcommand{\diagspace}{\mbox{\hspace{2em}}}
\newcommand{\ucontents}[2]{\addcontentsline{toc}{#1}{\numberline{}{#2}}}
\newcommand{\mb}[1]{\mathbf{#1}}
\newcommand{\mc}[1]{\mathcal{#1}}
\newcommand{\mr}[1]{\mathrm{#1}}
\newcommand{\cat}[1]{\mc{#1}}
\newcommand{\fcat}[1]{\mb{#1}}
\newcommand{\twid}[1]{\widetilde{#1}}
\newcommand{\slsh}{/\linebreak[0]}
\newcommand{\dt}{.\linebreak[0]}
\newcommand{\such}{\:|\:}
\newcommand{\epsln}{\varepsilon}
\newcommand{\implies}{\,\Rightarrow\,}
\newcommand{\op}{\mr{op}}
\newcommand{\Cat}{\fcat{Cat}}
\newcommand{\One}{\fcat{1}}
\newcommand{\Set}{\fcat{Set}}
\newcommand{\integers}{\mathbb{Z}}
\newcommand{\ftrcat}[2]{[#1,#2]}

\newcommand{\ogby}[1]{\lTo^{#1}}
\newcommand{\oppair}[2]{\pile{\rTo^{\scriptstyle #1}\\ 
\lTo_{\scriptstyle #2}}} 
\newcommand{\parpair}[2]{\pile{\rTo^{\scriptstyle #1}\\ 
\rTo_{\scriptstyle #2}}}
\newcommand{\parpairu}{\pile{\rTo\\ \rTo}}
\newarrow{Get}....>
\newarrow{Monic}>--->
\newcommand{\oppairu}{\pile{\rTo\\ \lTo}}
\newcommand{\reals}{\mathbb{R}}
\newcommand{\demph}[1]{\textbf{\textup{#1}}}
\newcommand{\done}{\hfill\ensuremath{\Box}}
\newenvironment{prooflike}[1]{\begin{trivlist}\item\textbf{#1}\ }
{\end{trivlist}}
\newenvironment{proof}{\begin{prooflike}{Proof}}{\end{prooflike}}
\newcommand{\scat}[1]{\mathbb{#1}}
\newcommand{\go}{\rTo\linebreak[0]}
\newcommand{\goby}[1]{\rTo^{#1}\linebreak[0]}
\newcommand{\iso}{\cong}
   
\newcommand{\eqv}{\simeq}
\newcommand{\url}[1]{#1}
\newcommand{\elt}[1]{\scat{E}\left(#1\right)}
\newcommand{\ob}{\mr{ob}\,}             
\newcommand{\of}{\,\raisebox{0.08ex}{\ensuremath{\scriptstyle\circ}}\,}
\newcommand{\sub}{\subseteq}
\newcommand{\ladj}{\dashv}     
\newcommand{\cell}[4]{\put(#1,#2){\makebox(0,0)[#3]{\ensuremath{#4}}}}
\newcommand{\rationals}{\mathbb{Q}}
\newcommand{\inc}[1]{R(#1)}
\newcommand{\incinf}[1]{\hat{R}(#1)}
\newcommand{\Dinj}[1]{\scat{D}^{\mr{inj}}_{#1}}
\newcommand{\binom}[2]{{#1 \choose #2}}
\newcommand{\bl}{{\scriptscriptstyle\bullet}}
\newcommand{\colimit}[1]{{\displaystyle\lim_{\longrightarrow}}\,}
\newcommand{\im}{\mr{im}}
\newcommand{\Aut}{\mr{Aut}}
\newarrow{Epic}----{>>}
\newcommand{\wkqt}{/\!/}
\newcommand{\Lef}[1]{L(#1)}
\newcommand{\Fix}[1]{\mb{Fix}\,#1}
\newcommand{\Alg}[1]{\mb{Alg}\,#1}
\newcommand{\Coalg}[1]{\mb{Coalg}\,#1}
\newcommand{\downset}[1]{\downarrow\! #1}

\newtheorem{thm}{Theorem}[section]
\newtheorem{propn}[thm]{Proposition}
\newtheorem{lemma}[thm]{Lemma}
\newtheorem{cor}[thm]{Corollary}
\newtheorem{predefn}[thm]{Definition}
\newenvironment{defn}{\begin{predefn}\upshape}{\end{predefn}}
\newtheorem{preexample}[thm]{Example}
\newenvironment{example}{\begin{preexample}\upshape}{\end{preexample}}
\newtheorem{preexamples}[thm]{Examples}
\newenvironment{examples}{\begin{preexamples}\upshape}{\end{preexamples}}

\begin{document}

\sloppy

\title{The Euler characteristic of a category}
\author{Tom Leinster%
\thanks{Department of Mathematics, University of Glasgow; 
tl@maths\protect\dt gla\protect\dt ac\protect\dt uk.  
Partially supported by a Nuffield Foundation award NUF-NAL 04 and an EPSRC
Advanced Research Fellowship.}
}
\date{}

\maketitle

\begin{center}
\textbf{Abstract}
\end{center}
The Euler characteristic of a finite category is defined and shown to be
compatible with Euler characteristics of other types of object, including
orbifolds.  A formula for the cardinality of the colimit of a diagram of sets
is proved, generalizing the classical inclusion-exclusion formula.  Both rest
on a generalization of M\"obius--Rota inversion from posets to categories.

\tableofcontents

\section*{Introduction}
\ucontents{section}{Introduction}
\label{sec:intro}

We first learn of Euler characteristic as `vertices minus edges plus faces',
and later as an alternating sum of ranks of homology groups.  But Euler
characteristic is much more fundamental than these definitions make apparent,
as has been made increasingly explicit over the last fifty years; it is
something akin to cardinality or measure.  More precisely, it is the
fundamental dimensionless quantity associated with an object.

Finite sets provide the very simplest context for Euler characteristic, and of
course the fundamental way to assign a quantity to a finite set is to count
its elements.  Indeed, Euler characteristic of topological spaces can usefully
be thought of as a generalization of cardinality; for instance, it obeys the
same laws with respect to unions and products.

A further example reinforces the point.  A subset of $\reals^n$ is
\demph{polyconvex} if it is a finite union of compact convex subsets.  Let
$V_n$ be the vector space of finitely additive measures, invariant under
Euclidean transformations, defined on the polyconvex subsets of $\reals^n$.
Hadwiger's Theorem~\cite{KR} states that $\dim V_n = n + 1$.  (See
also~\cite{SchWLP}.)  A natural basis consists of one $d$-dimensional measure
for each $d \in \{ 0, \ldots, n \}$: for instance, $\{ \textrm{Euler
characteristic}, \textrm{perimeter}, \textrm{area} \}$ when $n = 2$.  Thus, up
to scalar multiplication, Euler characteristic is the unique dimensionless
measure on polyconvex sets.

Schanuel~\cite{SchNSH} showed that in other contexts, Euler characteristic can
be defined in a way that makes its fundamental nature transparent.  He proved
that for a certain category of polyhedra, Euler characteristic is determined
by a straightforward universal property.

All of this makes clear the importance of defining and understanding Euler
characteristic in new contexts.  Here we do this for finite categories.  

One might envisage simply transporting the definition from spaces to
categories via the classifying space functor, as with other topological
invariants: that is, given a category $\scat{A}$, define $\chi(\scat{A})$ as
the Euler characteristic of the classifying space $B\scat{A}$.  The trouble
with this is that the Euler characteristic of $B\scat{A}$ is not always
defined.  Below we give a definition of the Euler characteristic of a category
that agrees with the topological Euler characteristic when the latter exists,
but is also valid in a range of situations when it does not.  It is a rational
number, not necessarily an integer.

A version of the definition can be given very succinctly.  Let $\scat{A}$ be a
finite category; totally order its objects as $a_1, \ldots, a_n$.  Let $Z$ be
the matrix whose $(i, j)$-entry is the number of arrows from $a_i$ to $a_j$.
Let $M = Z^{-1}$, assuming that this inverse exists.  Then $\chi(\scat{A})$ is
the sum of the entries of $M$.  Of course, it remains to convince the reader
that this definition is the right one.

The foundation on which this work rests is a generalization of
M\"obius--Rota inversion~(\S\ref{sec:mobius}).  Rota developed M\"obius
inversion for posets~\cite{Rota}; we develop it for categories.  (A poset is
viewed throughout as a category in which each hom-set has at most one element:
the objects are the elements of the poset, and there is an arrow $a \go b$ if
and only if $a \leq b$.)  This leads, among other things, to a `representation
formula': given any functor known to be a sum of representables, the formula
tells us the representation explicitly.  This in turn can be used to solve
enumeration problems, in the spirit of Rota's paper.

However, the main application of this generalized M\"obius inversion is to the
theory of the Euler characteristic of a category~(\S\ref{sec:euler}).  We
actually use a different definition than the one just given; it agrees with
the one above when $Z$ is invertible, but is valid for a wider class of
categories.  It depends on the idea of the `weight' of an object of a
category.  We justify the definition in two ways: by showing that it enjoys
the properties that the name would lead one to expect (behaviour with respect
to products, fibrations, etc.), and by demonstrating its compatibility with
Euler characteristics of other types of structure (groupoids, graphs,
topological spaces, orbifolds).

The technology of M\"obius inversion and weights also solves another problem:
what is the cardinality of a colimit?  For example, the union of a family of
sets and the quotient of a set by a free action of a group are both examples
of colimits of set-valued functors, and there are simple formulas for their
cardinalities.  (In the first case it is the inclusion-exclusion formula.)  We
generalize, giving a formula valid for any shape of
colimit~(\S\ref{sec:colimit}).

Rota and his school proved a large number of results on M\"obius inversion for
posets.  As we will repeatedly see, many are not truly order-theoretic: they
are facts about categories in general.  In particular, important theorems in
Rota's original work~\cite{Rota} generalize from posets to
categories~(\S\ref{sec:other}).

(The body of work on M\"obius inversion in finite lattices is not, however, so
ripe for generalization: a poset is a lattice just when the corresponding
category has products, but a finite category cannot have products unless it
is, in fact, a lattice.)

Other authors have considered different notions of M\"obius inversion for
categories; notably, there is that developed by Content, Lemay and
Leroux~\cite{CLL} and independently by Haigh~\cite{Hai}.  This generalizes
both Rota's notion for posets and Cartier and Foata's for monoids~\cite{CF}.
(Here a monoid is viewed as a one-object category.)  The relation between
their approach and ours is discussed in~\S\ref{sec:other}.  A third approach,
not discussed here, was taken by D\"ur~\cite{Dur}.

In the case of groupoids, our Euler characteristic of categories agrees with
Baez and Dolan's groupoid cardinality~\cite{BD}.  That in turn interacts well
with the species of Joyal~\cite{Joy,BLL}.  Our definition of the Lefschetz
number of an endofunctor~(\S\ref{sec:euler}) may perhaps be related to
Par\'e's definition of the cardinality of an endofunctor of the category of
finite sets~\cite{Pare}.

The view of Euler characteristic as generalized cardinality is promoted
in~\cite{SchNSH}, \cite{BD} and~\cite{ProEMG}.  The appearance of a
non-integral Euler characteristic is nothing new: see for instance
Wall~\cite{Wall}, Bass~\cite{Bass} and Cohen~\cite{Coh}, and the discussion of
orbifolds in~\S\ref{sec:euler}.

Ultimately it would be desirable to have the Euler characteristic of
categories described by a universal property, as Schanuel did for
polyhedra~\cite{SchNSH}.  For this, it may be necessary to relax the
constraints of the present work, where for simplicity our categories are
required to be finite and the coefficients are required to lie in the rig
(semiring) of rational numbers.  Rather than asking, as below, `does this
category have Euler characteristic (in $\rationals$)?', we should perhaps ask
`in what rig does the Euler characteristic of this category lie?'  However,
this is not pursued here.

\paragraph*{Acknowledgements} 
I thank John Baez, Andy Baker, Nick Gurski, Ieke Moerdijk and
Ivan Smith for inspiration and useful discussions.

\section{M\"obius inversion}
\label{sec:mobius}

We consider a finite category $\scat{A}$, writing $\ob\scat{A}$ for its set
of objects and, when $a$ and $b$ are objects, $\scat{A}(a, b)$ for the set of
maps from $a$ to $b$.

\begin{defn}    \label{defn:zeta-mobius}
We denote by $\inc{\scat{A}}$ the $\rationals$-algebra of functions
$\ob\scat{A} \times \ob\scat{A} \go \rationals$, with pointwise addition and
scalar multiplication, multiplication defined by
\[
(\theta \phi)(a, c)
=
\sum_{b \in \scat{A}}
\theta(a, b) \phi(b, c)
\]
($\theta, \phi \in \inc{\scat{A}}$, $a, c \in \scat{A}$), and the Kronecker
$\delta$ as unit.

The \demph{zeta function} $\zeta_{\scat{A}} = \zeta \in \inc{\scat{A}}$ is
defined by $\zeta(a, b) = |\scat{A}(a, b)|$.  If $\zeta$ is invertible in
$\inc{\scat{A}}$ then $\scat{A}$ is said to \demph{have M\"obius
inversion}; its inverse $\mu_\scat{A} = \mu = \zeta^{-1}$ is the
\demph{M\"obius function} of $\scat{A}$.
\end{defn}

If a total ordering is chosen on the $n$ objects of $\scat{A}$ then
$\inc{\scat{A}}$ can be regarded as the algebra of $n\times n$ matrices over
$\rationals$.  The defining equations of the M\"obius function are
\[
\sum_b \mu(a, b) \zeta(b, c)
=
\delta(a, c)
=
\sum_b \zeta(a, b) \mu(b, c)
\]
for all $a, c \in \scat{A}$.  By finite-dimensionality, either one implies
the other.

The definitions above could be made for directed graphs rather than
categories, since they do not refer to composition.  However, this generality
seems to be inappropriate.  For example, the definition of M\"obius inversion
will lead to a definition of Euler characteristic, and if we use graphs rather
than categories then we obtain something other than `vertices minus edges'.
Proposition~\ref{propn:euler-graph} clarifies this point.

A different notion of M\"obius inversion for categories has been considered;
see~\S\ref{sec:other}.

\begin{examples}    \label{egs:mobius}
\begin{enumerate}
\item   \label{eg:mobius-posets}
Any finite poset $A$ has M\"obius inversion; this special case was
investigated by Rota~\cite{Rota} and others.  We may compute $\mu(a, c)$ by
induction on the number of elements between $a$ and $c$:
\[
\mu(a, c)
=
\delta(a, c) - \sum_{b: a \leq b < c} \mu(a, b).
\]
In particular, $\mu(a, c) = 0$ unless $a \leq c$, and $\mu(a, a) = 1$ for
all $a$.

\item Let $M$ be a finite monoid, regarded as a category with unique object
$\star$.  Then $\zeta(\star, \star) = |M|$, so $\mu(\star, \star) = 1/|M|$.

\item   \label{eg:mobius-Dinj}
Let $N \geq 0$.  Write $\Dinj{N}$ for the category with objects $0, \ldots, N$
whose maps $a \go b$ are the order-preserving injections $\{1, \ldots, a\} \go
\{1, \ldots, b\}$.  Then $\zeta(a, b) = \binom{b}{a}$, and it is easily
checked that $\mu(a, b) = (-1)^{b - a} \binom{b}{a}$.  If we use surjections
instead of injections then $\zeta(a, b) = \binom{a - 1}{b - 1}$ and $\mu(a, b)
= (-1)^{a - b} \binom{a - 1}{b - 1}$.  
\end{enumerate}
\end{examples}

A category with M\"obius inversion must be skeletal, for otherwise the
matrix of $\zeta$ would have two identical rows.  The property of having
M\"obius inversion is not, therefore, invariant under equivalence of
categories.

In general we cannot hope to just spot the M\"obius function of a category.
In~\ref{lemma:idem-id}--\ref{eg:finsets} we build tools for computing M\"obius
functions.  These cover large classes of categories, although not every finite
skeletal category has M\"obius
inversion~(\ref{egs:wtgs}(\ref{eg:wtgs-none}),~(\ref{eg:wtgs-many})).

Let $n \geq 0$, let $\scat{A}$ be a category or a directed graph, and let
$a, b \in \scat{A}$. An \demph{$n$-path} from $a$ to $b$ is a diagram
\begin{equation}    \label{eq:path}
a = a_0 \goby{f_1} a_1 \goby{f_2}
\ \cdots \
\goby{f_n} a_n = b
\end{equation}
in $\scat{A}$.  It is a \demph{circuit} if $a = b$, and (when $\scat{A}$ is a
category) \demph{nondegenerate} if no $f_i$ is an identity.

\begin{lemma}   \label{lemma:idem-id}
The following conditions on a finite category $\scat{A}$ are equivalent:
\begin{enumerate}
\item   \label{item:idem-id}
every idempotent in $\scat{A}$ is an identity
\item   \label{item:endo-auto}
every endomorphism in $\scat{A}$ is an automorphism
\item   \label{item:circuit-free}
every circuit in $\scat{A}$ consists entirely of isomorphisms.
\end{enumerate}
\end{lemma}

\begin{proof}
(\ref{item:idem-id}) $\implies$ (\ref{item:endo-auto}) follows from the fact
that if $f$ is an element of a finite monoid then some positive power of $f$
is idempotent.  The other implications are straightforward.  \done
\end{proof}

\begin{thm}     \label{thm:mobius-idem-id}
Let $\scat{A}$ be a finite skeletal category in which the only idempotents
are identities.  Then $\scat{A}$ has M\"obius inversion given by
\[
\mu(a, b)
=
\sum (-1)^n/|\Aut(a_0)| \cdots |\Aut(a_n)|
\]
where $\Aut(a)$ is the automorphism group of $a \in \scat{A}$ and the sum
runs over all $n \geq 0$ and paths~(\ref{eq:path}) for which $a_0, \ldots,
a_n$ are all distinct.
\end{thm}

\begin{proof}
First observe that for a path~(\ref{eq:path}) in $\scat{A}$, if $a_0 \neq
a_1 \neq \cdots \neq a_n$ then the $a_i$s are all distinct.  Indeed,
if $0 \leq i < j \leq n$ and $a_i = a_j$ then the sub-path running from
$a_i$ to $a_j$ is a circuit, so by Lemma~\ref{lemma:idem-id}, $f_{i + 1}$
is an isomorphism, and by skeletality, $a_i = a_{i + 1}$.

Now let $a, c \in \scat{A}$ and define $\mu$ by the formula above.  We have
\begin{eqnarray*}
\sum_{b \in \scat{A}} \mu(a, b) \zeta(b, c) &=  &
\mu(a, c) \zeta(c, c)
+
\sum_{b: b \neq c} \mu(a, b) \zeta(b, c)    \\
    &=  &
|\Aut(c)|
\left\{
\mu(a, c)
+
\sum_{b: b \neq c,
\
g \in \scat{A}(b, c)}
\mu(a, b) 
/ |\Aut(c)|
\right\},
\end{eqnarray*}
and by definition of $\mu$, the term in braces collapses to $\delta(a,
c)/|\Aut(a)|$.
\done
\end{proof}

\begin{cor} \label{cor:mobius-circuit-free}
Let $\scat{A}$ be a finite skeletal category in which the only
endomorphisms are identities.  Then $\scat{A}$ has M\"obius inversion given
by
\[
\mu(a, b) = \sum_{n \geq 0} (-1)^n | \{ \textrm{nondegenerate 
$n$-paths from $a$ to $b$} \} |
\in \integers.
\]
\done
\end{cor}

When $\scat{A}$ is a poset, this is Philip Hall's theorem
(Proposition~3.8.5 of~\cite{Sta} and Proposition~6 of~\cite{Rota}). 

Recall that an \demph{epi-mono factorization system} $(\cat{E}, \cat{M})$ on a
category $\scat{A}$ consists of a class $\cat{E}$ of epimorphisms in
$\scat{A}$ and a class $\cat{M}$ of monomorphisms in $\scat{A}$, satisfying
axioms~\cite{FK}.  The axioms imply that every map $f$ in $\scat{A}$ can be
expressed as $me$ for some $e \in \cat{E}$ and $m \in \cat{M}$, and that this
factorization is essentially unique: the other pairs $(e', m') \in \cat{E}
\times \cat{M}$ satisfying $m'e' = f$ are those of the form $(ie, m i^{-1})$
where $i$ is an isomorphism.

\begin{thm}     \label{thm:mobius-fs}
Let $\scat{A}$ be a finite skeletal category with an epi-mono
factorization system $(\cat{E}, \cat{M})$.  Then $\scat{A}$ has M\"obius
inversion given by
\[
\mu(a, b) 
= 
\sum (-1)^n / |\Aut(a_0)| \cdots |\Aut(a_n)|
\]
where the sum is over all $n \geq r \geq 0$ and paths~(\ref{eq:path}) such
that $a_0, \ldots, a_r$ are distinct, $a_r, \ldots, a_n$ are distinct,
$f_1, \ldots, f_r \in \cat{M}$, and $f_{r + 1}, \ldots, f_n \in \cat{E}$.
\end{thm}

\begin{proof}
The objects of $\scat{A}$ and the arrows in $\cat{E}$ determine a
subcategory of $\scat{A}$, also denoted $\cat{E}$; it satisfies the
hypotheses of Theorem~\ref{thm:mobius-idem-id} and therefore has M\"obius
inversion.  The same is true of $\cat{M}$.

Any element $\alpha \in \rationals^{\ob\scat{A}}$ gives rise to an element
of $\inc{\scat{A}}$, also denoted $\alpha$ and defined by $\alpha(a, b) =
\delta(a, b) \alpha(b)$.  The resulting map from $\rationals^{\ob\scat{A}}$
to $\inc{\scat{A}}$ preserves multiplication (where the multiplication on
$\rationals^{\ob\scat{A}}$ is pointwise).  We have elements $|\Aut|$ and
$1/|\Aut|$ of $\rationals^{\ob\scat{A}}$, where, for instance, $|\Aut|(a) =
|\Aut(a)|$.

By the essentially unique factorization property,
$\zeta_{\scat{A}} = \zeta_\cat{E} \cdot \frac{1}{|\Aut|} \cdot
\zeta_\cat{M}$.  Hence $\scat{A}$ has M\"obius function
$\mu_{\scat{A}} = \mu_{\cat{M}} \cdot |\Aut| \cdot
\mu_{\cat{E}}$.  Theorem~\ref{thm:mobius-idem-id} then gives the
formula claimed.  \done
\end{proof}

\begin{example} \label{eg:finsets}
Let $N \geq 0$ and write $\scat{F}_N$ for the full subcategory of $\Set$ with
objects $1, \ldots, N$, where $n$ denotes a (chosen) $n$-element set.  Let
$\cat{E}$ be the set of surjections in $\scat{F}_N$ and $\cat{M}$ the set of
injections; then $(\cat{E}, \cat{M})$ is a factorization system.
Theorem~\ref{thm:mobius-fs} gives a formula for the inverse of the matrix
$(i^j)_{i, j}$.  For instance, take $N = 3$; then $\mu(1, 2)$ may be computed
as follows:
\[
\begin{tabular}{ll}
Paths                                   &
Contribution to sum                     \\
$1 \rMonic^2 2$                         &
$-2/1!2! = -1$                          \\
$1 \rMonic^3 3 \rEpic^6 2$              &
$3 \cdot 6/1!3!2! = 3/2$                \\
$1 \rMonic^2 2 \rMonic^6 3 \rEpic^6 2$  &
$-2 \cdot 6 \cdot 6/1!2!3!2! = -3$
\end{tabular}
\]
Here `$1 \rMonic^2 2$' means that there are $2$ monomorphisms from $1$ to
$2$, `$3 \rEpic^6 2$' that there are $6$ epimorphisms from $3$ to $2$, etc.
Hence $\mu(1, 2) = -1 + 3/2 - 3 = -5/2$.
\end{example}

One of the uses of the M\"obius function is to calculate Euler
characteristic~(\S\ref{sec:euler}).  Another is to calculate representations.
Specifically, suppose that we have a $\Set$-valued functor known to be
\demph{familially representable}, that is, a coproduct of representables.  The
Yoneda Lemma tells us that the family of representing objects is unique (up to
isomorphism); better, if we have M\"obius inversion, there is actually a
formula for it:
\begin{propn}   \label{propn:fr-coeffs}
Let $\scat{A}$ be a finite category with M\"obius inversion and let $X:
\scat{A} \go \Set$ be a functor satisfying
\[
X 
\iso
\sum_a r(a) \scat{A}(a, \dashbk)
\]
for some natural numbers $r(a)$ ($a \in \scat{A}$).  Then 
\[
r(a) = \sum_b |Xb| \mu(b, a)
\]
for all $a \in \scat{A}$.
\end{propn}

\begin{proof}
Follows from the definition of M\"obius function. 
\done
\end{proof}

In the spirit of Rota's programme, this can be applied to solve counting
problems, as illustrated by the following standard example.

\begin{example}
A \demph{derangement} is a permutation without fixed points.  We calculate
$d_n$, the number of derangements of $n$ letters.

Fix $N \geq 0$.  Take the category $\Dinj{N}$ of
Example~\ref{egs:mobius}(\ref{eg:mobius-Dinj}) and the functor $S: \Dinj{N}
\go \Set$ defined as follows: $S(n)$ is $S_n$, the underlying set of the
$n$th symmetric group, and if $f \in \Dinj{N}(m, n)$ and $\tau \in S_m$,
the induced permutation $S_f(\tau) \in S_n$ acts as $\tau$ on $\im f$ and
fixes all other points.  Any permutation consists of a derangement together
with some fixed points, so
\[
S_n
\iso 
\sum_m d_m \Dinj{N}(m, n).
\]
Then by Proposition~\ref{propn:fr-coeffs},
\[
d_n
=
\sum_m |S_m| \mu(m, n)
=
\sum_{m} m! (-1)^{n - m} \binom{n}{m}
=
n! \left( \frac{1}{0!} - \frac{1}{1!} + \cdots + \frac{(-1)^n}{n!} \right).
\]
\end{example}

To set up the theory of Euler characteristic we will not need the full
strength of M\"obius invertibility; the following suffices.

\begin{defn}
Let $\scat{A}$ be a finite category.  A \demph{weighting} on $\scat{A}$ is
a function $k^\bl: \ob\scat{A} \go \rationals$ such that for all $a \in
\scat{A}$,
\[
\sum_b \zeta(a, b) k^b = 1.
\]
A \demph{coweighting} $k_\bl$ on $\scat{A}$ is a weighting on
$\scat{A}^\op$.
\end{defn}

Note that $\scat{A}$ has M\"obius inversion if and only if it has a unique
weighting, if and only if it has a unique coweighting; they are given by
\[
k^a
=
\sum_b \mu(a, b),
\qquad
k_b
=
\sum_a \mu(a, b).
\]

\begin{examples}    \label{egs:wtgs}
\begin{enumerate}
\item   \label{eg:wtgs-pushout}
Let $\scat{L}$ be the category
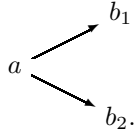
\[
\begin{diagram}[height=1em]
        &       &b_1    \\
        &\ruTo  &       \\
a       &       &       \\
        &\rdTo  &       \\
        &       &b_2.   \\
\end{diagram}
\]
Then the unique weighting $k^\bl$ on $\scat{L}$ is $(k^a, k^{b_1}, k^{b_2})
= (-1, 1, 1)$. 

\item   \label{eg:wtgs-monoid}
Let $M$ be a finite monoid, regarded as a category with unique object
$\star$.  Again there is a unique weighting $k^\bl$, with $k^\star =
1/|M|$.

\item   \label{eg:wtgs-terminal}
If $\scat{A}$ has a terminal object $1$ then $\delta(\dashbk, 1)$ is a 
weighting on $\scat{A}$.

\item   \label{eg:wtgs-none}
A finite category may admit no weighting at all (even if Cauchy-complete).
An example is the category $\scat{A}$ with objects and arrows
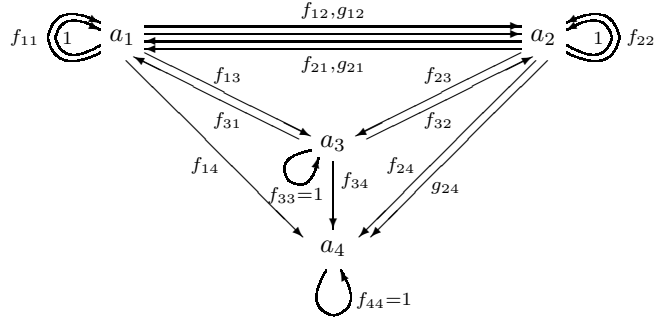
\[
\setlength{\unitlength}{1mm}
\begin{picture}(80,42)(-40,-9)
\cell{-28}{28}{c}{a_1}
\cell{28}{28}{c}{a_2}
\cell{0}{14}{c}{a_3}
\cell{0}{0}{c}{a_4}
\put(-25,29.5){\vector(1,0){50}}        
\put(-25,28.5){\vector(1,0){50}}        
\put(25,27.5){\vector(-1,0){50}}        
\put(25,26.5){\vector(-1,0){50}}        
\put(-25,26){\vector(2,-1){22}}         
\put(-4.5,14.5){\vector(-2,1){22}}      
\put(25,26){\vector(-2,-1){22}}         
\put(4.5,14.5){\vector(2,1){22}}        
\put(0,11.5){\vector(0,-1){9}}          
\put(-27.5,25){\vector(1,-1){23.5}}     
\put(27,25){\vector(-1,-1){23.5}}     
\put(28.5,25){\vector(-1,-1){23.5}}     
\qbezier(-31,27)(-36,24)(-37,28)
\qbezier(-31,29)(-36,32)(-37,28)
\put(-31,29){\vector(2,-1){0}}          
\qbezier(-31,26)(-37,23)(-38,28)
\qbezier(-31,30)(-37,33)(-38,28)
\put(-31,30){\vector(2,-1){0}}          
\qbezier(31,27)(36,24)(37,28)
\qbezier(31,29)(36,32)(37,28)
\put(31,29){\vector(-2,-1){0}}          
\qbezier(31,26)(37,23)(38,28)
\qbezier(31,30)(37,33)(38,28)
\put(31,30){\vector(-2,-1){0}}          
\qbezier(-1,-3)(-4,-8)(0,-9)
\qbezier(1,-3)(4,-8)(0,-9)
\put(1,-3){\vector(-1,2){0}}            
\qbezier(-3,13)(-8,12)(-6,9)
\qbezier(-2,12)(-3,7)(-6,9)
\put(-2,12){\vector(1,4){0}}            
\cell{0}{30.5}{b}{\scriptstyle{f_{12}, g_{12}}}
\cell{0}{25.5}{t}{\scriptstyle{f_{21}, g_{21}}}
\cell{-14}{23}{c}{\scriptstyle{f_{13}}}
\cell{-14}{17}{c}{\scriptstyle{f_{31}}}
\cell{14}{23}{c}{\scriptstyle{f_{23}}}
\cell{14}{17}{c}{\scriptstyle{f_{32}}}
\cell{1}{9}{l}{\scriptstyle{f_{34}}}
\cell{-17}{11}{c}{\scriptstyle{f_{14}}}
\cell{9}{11}{c}{\scriptstyle{f_{24}}}
\cell{15}{8}{c}{\scriptstyle{g_{24}}}
\cell{-36}{28}{l}{\scriptstyle{1}}
\cell{-39}{28}{r}{\scriptstyle{f_{11}}}
\cell{36}{28}{r}{\scriptstyle{1}}
\cell{39}{28}{l}{\scriptstyle{f_{22}}}
\cell{3}{-7}{l}{\scriptstyle{f_{44} = 1}}
\cell{-5}{7}{c}{\scriptstyle{f_{33} = 1}}
\end{picture}
\]
where if $a_i \goby{p} a_j \goby{q} a_k$ and
neither $p$ nor $q$ is an identity then $q \of p = f_{ik}$. 

\item   \label{eg:wtgs-many}
A category may certainly have more than one weighting: for instance, if
$\scat{A}$ is the category consisting of two objects and a single isomorphism
between them, a weighting on $\scat{A}$ is any pair of rational numbers whose
sum is $1$.  But even a skeletal category may admit more than one weighting.
Indeed, the full subcategories $\scat{B} = \{a_1, a_2\}$ and $\scat{C} =
\{a_1, a_2, a_3\}$ of the category $\scat{A}$ of the previous example both
have infinitely many weightings.
\end{enumerate}
\end{examples}

In contrast to M\"obius invertibility, the property of admitting at least
one weighting is invariant under equivalence:

\begin{lemma}   \label{lemma:eqv-wtg}
Let $\scat{A}$ and $\scat{B}$ be equivalent finite categories.  Then
$\scat{A}$ admits a weighting if and only if $\scat{B}$ does.
\end{lemma}

\begin{proof}
Let $F: \scat{A} \go \scat{B}$ be an equivalence.  Given $a \in \scat{A}$,
write $C_a$ for the number of objects in the isomorphism class of $a$, and
similarly $C_b$ when $b \in \scat{B}$.
Take a weighting $l^\bl$ on $\scat{B}$ and put $k^a = (C_{Fa}/C_a)
l^{Fa}$; then $k^\bl$ is a weighting on $\scat{A}$.
\done
\end{proof}

Weightings and M\"obius functions are compatible with sums and products of
categories; the following lemma is easily verified.
\begin{lemma}   \label{lemma:mobius-sum-prod}
Let $n \geq 0$ and let $\scat{A}_1, \ldots, \scat{A}_n$ be finite
categories.
\begin{enumerate}
\item If each $\scat{A}_i$ has a weighting $k^\bl_i$ then $\sum_i
\scat{A}_i$ has a weighting $l^\bl$ given by $l^a = k^a_i$  whenever $a
\in \scat{A}_i$.  If each $\scat{A}_i$ has M\"obius inversion then so does
$\sum_i \scat{A}_i$, where for $a \in \scat{A}_i$ and $b \in \scat{A}_j$,
\[
\mu(a, b)
=
\left\{
\begin{array}{ll}
\mu_{\scat{A}_i} (a, b) &\textrm{if } i = j     \\
0                       &\textrm{otherwise}.
\end{array}
\right.
\]
\item If each $\scat{A}_i$ has a weighting $k^\bl_i$ then $\prod_i
\scat{A}_i$ has a weighting $l^\bl$ given by
$
l^{(a_1, \ldots, a_n)}
=
k_1^{a_1} \cdots k_n^{a_n}.
$
If each $\scat{A}_i$ has M\"obius inversion then so does $\prod_i
\scat{A}_i$, with
\[
\mu((a_1, \ldots, a_n), (b_1, \ldots, b_n))
=
\mu(a_1, b_1) \cdots \mu(a_n, b_n).
\]
\done
\end{enumerate}
\end{lemma}

Recall that any functor $X$ taking values in $\Set$ or $\Cat$ has a category
of elements $\elt{X}$, and that in the $\Cat$-valued case, this applies
equally when $X$ is a weak (or `pseudo') functor; see the Appendix.  We call
$X$ \demph{finite} if $\elt{X}$ is finite.  When the domain category
$\scat{A}$ is finite, this just means that each set or category $Xa$ is
finite.

\begin{lemma}   \label{lemma:wtg-fib}
Let $\scat{A}$ be a finite category and $X: \scat{A} \go \Cat$ a finite
weak functor.  Suppose that we have weightings on $\scat{A}$ and on each $Xa$,
all written $k^\bl$.  Then there is a weighting on $\elt{X}$ defined by
$k^{(a, x)} = k^a k^x$ ($a \in \scat{A}$, $x \in Xa$).
\end{lemma}

\begin{proof}
Let $a \in \scat{A}$ and $x \in Xa$.  Then
\begin{eqnarray*}
\sum_{(b, y) \in \elt{X}} \zeta((a, x), (b, y)) k^b k^y &
=       &
\sum_b \sum_{f \in \scat{A}(a, b)}
\left( \sum_{y \in Xb} \zeta((Xf)x, y) k^y \right) k^b  \\
        &=      &
\sum_b \zeta(a, b) k^b  
= 
1.
\end{eqnarray*}
\done
\end{proof}

This result will be used to show how Euler characteristic behaves
with respect to fibrations.

\section{Euler characteristic}
\label{sec:euler}

In this section, the Euler characteristic of a category is defined and its
basic properties are established.  The definition is justified by a series of
propositions showing its compatibility with the Euler characteristics of other
types of object: graphs, topological spaces, and orbifolds.  There follows
a brief discussion of the Lefschetz number of an endofunctor.

\begin{lemma}   \label{lemma:wtg-cowtg}
Let $\scat{A}$ be a finite category, $k^\bl$ a weighting on
$\scat{A}$, and $k_\bl$ a coweighting on $\scat{A}$.  Then
$\sum_a k^a = \sum_a k_a$.
\end{lemma}

\begin{proof}
\[
\sum_b k^b
=
\sum_b \left( \sum_a k_a \zeta(a, b) \right) k^b
=
\sum_a k_a \left( \sum_b \zeta(a, b) k^b \right)
=
\sum_a k_a.
\]
\done
\end{proof}

\begin{defn}
A finite category $\scat{A}$ \demph{has Euler characteristic} if it admits
both a weighting and a coweighting.  Its \demph{Euler characteristic} is
then
\[
\chi(\scat{A}) = \sum_a k^a = \sum_a k_a
\]
for any weighting $k^\bl$ and coweighting $k_\bl$.
\end{defn}

Any category $\scat{A}$ with M\"obius inversion has Euler characteristic,
$\chi(\scat{A}) = \sum_{a, b} \mu(a, b)$, as in the Introduction.  

\begin{examples}    \label{egs:euler-chars}
\begin{enumerate}
\item   \label{eg:euler-discrete}
If $\scat{A}$ is a finite discrete category then $\chi(\scat{A}) =
|\ob\scat{A}|$.

\item   \label{eg:euler-monoid}
If $M$ is a finite monoid then $\chi(M) = 1/|M|$.  When $M$ is a group, this
can be understood as follows: $M$ acts freely on the contractible space $EM$,
which has Euler characteristic $1$; one would therefore expect the quotient
space $BM$ to have Euler characteristic $1/|M|$.  (Compare~\cite{Wall}
and~\cite{Coh}.)

\item   \label{eg:euler-poset} 
By Corollary~\ref{cor:mobius-circuit-free}, a finite poset $A$ has Euler
characteristic $\sum_{n \geq 0} (-1)^n c_n$, where $c_n$ is the number of
chains in $A$ of length $n$.  (See \cite{Pul}, \cite{Fol}, \cite{Rota} and
\cite{Far} for connections with poset homology, and~\S\ref{sec:other} for
further comparisons with the Rota theory.)  More generally, the results
of~\S\ref{sec:mobius} give formulas for the Euler characteristic of any finite
category that either has no non-trivial idempotents or admits an epi-mono
factorization system.

\item \label{eg:euler-init-term} If $\scat{A}$ has Euler characteristic and
either an initial or a terminal object then $\chi(\scat{A}) = 1$,
by~\ref{egs:wtgs}(\ref{eg:wtgs-terminal}) and its dual.  In particular, any
finite category with both an initial and a terminal object has Euler
characteristic $1$.  This applies, for instance, to the category $\scat{C}$
of~\ref{egs:wtgs}(\ref{eg:wtgs-many}).  Hence having M\"obius inversion is a
strictly stronger property than having Euler characteristic.

\item Euler characteristic is not invariant under Morita equivalence.  For
example, the two-element monoid consisting of the identity and an idempotent
is Morita equivalent to the category freely generated by a split epimorphism,
but their respective Euler characteristics are~$1/2$ and~$1$.
\end{enumerate}
\end{examples}

Clearly $\chi(\scat{A}^\op) = \chi(\scat{A})$, one side being defined when
the other is.  The next few propositions set out further basic properties
of Euler characteristic.

\begin{propn}   \label{propn:euler-adjn-eqv}
Let $\scat{A}$ and $\scat{B}$ be finite categories.
\begin{enumerate}
\item   \label{item:euler-adjn}
If there is an adjunction $\scat{A} \oppairu \scat{B}$ and both
$\scat{A}$ and $\scat{B}$ have Euler characteristic then $\chi(\scat{A})
= \chi(\scat{B})$.
\item   \label{item:euler-eqv}
If $\scat{A} \eqv \scat{B}$ then $\scat{A}$ has Euler characteristic
  if and only if $\scat{B}$ does, and in that case $\chi(\scat{A}) =
  \chi(\scat{B})$.
\end{enumerate}
\end{propn}

In~(\ref{item:euler-adjn}), it may be that one category has Euler
characteristic but the other does not: consider, for instance, the unique
functor from the category $\scat{A}$ of~\ref{egs:wtgs}(\ref{eg:wtgs-none}) to
the terminal category.

\begin{proof}
\begin{enumerate}
\item Suppose that $\scat{A} \oppair{F}{G} \scat{B}$ with $F \ladj G$.
Then $\zeta(Fa, b) = \zeta(a, Gb)$ for all $a \in \scat{A}$, $b \in
\scat{B}$; write $\zeta(a, b)$ for their common value.  Take a coweighting
$k_\bl$ on $\scat{A}$ and a weighting $k^\bl$ on $\scat{B}$.  Then $\sum_a
k_a = \sum_b k^b$ by the same proof as that of Lemma~\ref{lemma:wtg-cowtg}.

\item The first statement follows from Lemma~\ref{lemma:eqv-wtg} and its
dual, and the second from~(\ref{item:euler-adjn}). 
\done
\end{enumerate}
\end{proof}

\begin{example} \label{eg:euler-i-t-functor}
If $\scat{B}$ is a category with an initial or a terminal object then
$\chi(\scat{A}^\scat{B}) = \chi(\scat{A})$ for all $\scat{A}$, provided that
both Euler characteristics exist.  Indeed, if $0$ is initial in $\scat{B}$
then evaluation at $0$ is right adjoint to the diagonal functor $\scat{A} \go
\scat{A}^\scat{B}$.
\end{example}

\begin{propn}   \label{propn:euler-sum-prod}
Let $n \geq 0$ and let $\scat{A}_1, \ldots, \scat{A}_n$ be finite
categories that all have Euler characteristic.  Then $\sum_i \scat{A}_i$
and $\prod_i \scat{A}_i$ have Euler characteristic, with
\[
\chi \left( \sum_i \scat{A}_i \right)
=
\sum_i \chi(\scat{A}_i),
\qquad
\chi \left( \prod_i \scat{A}_i \right)
=
\prod_i \chi(\scat{A}_i).
\]
\end{propn}

\begin{proof}
Follows from Lemma~\ref{lemma:mobius-sum-prod}.
\done
\end{proof}

\begin{example}       \label{eg:groupoid-chi}
Let $\scat{A}$ be a finite groupoid.  Choose one object $a_i$ from each
connected-component of $\scat{A}$, and write $G_i$ for the automorphism group
of $a_i$.  Then $\scat{A} \eqv \sum_i G_i$, so
by~\ref{egs:euler-chars}(\ref{eg:euler-monoid}),
\ref{propn:euler-adjn-eqv}(\ref{item:euler-eqv})
and~\ref{propn:euler-sum-prod}, we have $\chi(\scat{A}) = \sum_i 1/|G_i|$.
This is what Baez and Dolan call the cardinality of the groupoid
$\scat{A}$~\cite{BD}.
\end{example}

One might also ask whether $\chi(\scat{A}^\scat{B}) =
\chi(\scat{A})^{\chi(\scat{B})}$.  By
\ref{egs:euler-chars}(\ref{eg:euler-init-term}),
\ref{eg:euler-i-t-functor} and~\ref{propn:euler-sum-prod}, the
answer is yes if every connected-component of $\scat{B}$ has an
initial or a terminal object (and all the Euler characteristics
exist).  But in general the answer is no: for instance, take
$\scat{A}$ to be the $2$-object discrete category and $\scat{B}$
to be the category
of~\ref{egs:colim-formula}(\ref{eg:colim-coeq}).  See also
Propp~\cite{ProEEM}, Speed~\cite{Spe}, and \S5,~6 of
Rota~\cite{Rota}.

An important property of topological Euler characteristic is its behaviour
with respect to fibrations.  Let $E \go A$ be a topological fibration.  If $A$
has connected-components $A_1, \ldots, A_n$ and $X_i$ is the fibre in the
$i$th component then, under suitable hypotheses, $\chi(E) = \sum_i
\chi(A_i)\chi(X_i)$.  Lemma~\ref{lemma:wtg-fib} implies the following
categorical analogue.

\begin{propn}   \label{propn:euler-fib}
Let $\scat{A}$ be a finite category and $X: \scat{A} \go \Cat$ a finite
weak functor.  Let $k^\bl$ be a weighting on $\scat{A}$ and suppose that
$\elt{X}$ and each $Xa$ have Euler characteristic.  Then
\[
\chi(\elt{X})
=
\sum_a k^a \chi(Xa).
\]
\done
\end{propn}

\begin{examples}
\begin{enumerate}
\item When $X$ is a finite $\Set$-valued functor, $\chi(\elt{X}) = \sum_a k^a
|Xa|$.  For example, let $M$ be a finite monoid.  A finite functor $X: M \go
\Set$ is a finite set $S$ with a left $M$-action.  Following~\cite{BD}, we
write $\elt{X}$ as $S\wkqt M$, the \demph{weak quotient} of $S$ by $M$.  (Its
objects are the elements of $S$, and the arrows $s \go s'$ are
the elements $m \in M$ satisfying $ms = s'$.)  Then $\chi(S\wkqt M) =
|S|/|M|$.

\item
Define a sequence $(\scat{S}^n)_{n \geq -1}$ of categories inductively as
follows.  $\scat{S}^{-1}$ is empty.  Let $\scat{L}$ be the category
of~\ref{egs:wtgs}(\ref{eg:wtgs-pushout}); define $X: \scat{L} \go \Cat$ by
$Xa = \scat{S}^{n - 1}$ and $Xb_1 = Xb_2 = \One$ (the terminal category);
put $\scat{S}^n = \elt{X}$.  Explicitly, $\scat{S}^n$ is the poset
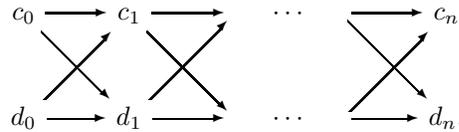
\[
\begin{diagram}
c_0     &\rTo   &c_1    &\rTo   &\mbox{}&\cdots &\mbox{}&\rTo   &c_n    \\
        &
\rdTo\ruTo      &       &
\rdTo\ruTo                      &       &       &       &
\rdTo\ruTo                                                      &       \\
d_0     &\rTo   &d_1    &\rTo   &\mbox{}&\cdots &\mbox{}&\rTo   &d_n.   \\
\end{diagram}
\]
(If we take the usual expression of the
topological $n$-sphere $S^n$ as a CW-complex with two cells in each
dimension $\leq n$ then $\scat{S}^n$ is the set of cells ordered by
inclusion; $S^n$ is the classifying space of $\scat{S}^n$.)

Each $\scat{S}^n$ is a poset, so has Euler characteristic.  By
Proposition~\ref{propn:euler-fib}, 
\[
\chi(\scat{S}^n) 
= 
-\chi(\scat{S}^{n -1}) + 2\chi(\One) 
= 
2 - \chi(\scat{S}^{n - 1})
\]
for all $n \geq 0$; also $\chi(\scat{S}^{-1}) = 0$.  Hence
$\chi(\scat{S}^n) = 1 + (-1)^n$.

\end{enumerate}
\end{examples}

The next three propositions show how the Euler characteristics of various
types of structure are compatible with that of categories.

First, Euler characteristic of categories extends Euler characteristic of
graphs.  More precisely, let $G = (G_1 \parpairu G_0)$ be a directed graph,
where $G_1$ is the set of edges and $G_0$ the set of vertices.  We will show
that if $F(G)$ is the free category on $G$ then $\chi(F(G)) = |G_0| - |G_1|$.
This only makes sense if $F(G)$ is finite, which is the case if and only if
$G$ is finite and circuit-free; then $F(G)$ is also circuit-free.  (A directed
graph is \demph{circuit-free} if it contains no circuits of non-zero length,
and a category is \demph{circuit-free} if every circuit consists entirely of
identities.)

\begin{propn}   \label{propn:euler-graph}
Let $G$ be a finite circuit-free directed graph.  Then $\chi(F(G))$ is
defined and equal to $|G_0| - |G_1|$.
\end{propn}

\begin{proof}
Given $a, b \in G_0$, write $\zeta_G(a, b)$ for the number of edges from
$a$ to $b$ in $G$.  Then $\zeta_{F(G)} = \sum_{n \geq 0} \zeta_G^n$ in
$\inc{F(G)}$, the sum being finite since $G$ is circuit-free.  Hence
$\mu_{F(G)} = \delta - \zeta_G$, and the result follows.
\done
\end{proof}

This suggests that in the present context, it is more fruitful to view a
graph as a special category (via $F$) than a category as a graph with
structure.  Compare the comments after Definition~\ref{defn:zeta-mobius}.

The second result compares the Euler characteristics of categories and
topological spaces.  We show that under suitable hypotheses, $\chi(B\scat{A})
= \chi(\scat{A})$, where $B\scat{A}$ is the classifying space of a category
$\scat{A}$ (that is, the geometric realization of its nerve $N\scat{A}$).  To
ensure that $B\scat{A}$ has Euler characteristic, we assume that $N\scat{A}$
contains only finitely many nondegenerate simplices; then
\[
\chi(B\scat{A})
=
\sum_{n \geq 0}
(-1)^n
| \{ \textrm{nondegenerate $n$-simplices in $N\scat{A}$} \} |.
\]
An $n$-simplex in $N\scat{A}$ is just an $n$-path in $\scat{A}$, and is
nondegenerate in the sense of simplicial sets if and only if it is
nondegenerate as a path, so $\scat{A}$ must contain only finitely many
nondegenerate paths.  This is the case if and only if $\scat{A}$ is
circuit-free, if and only if $\scat{A}$ is skeletal and contains no
endomorphisms except identities.  So by
Corollary~\ref{cor:mobius-circuit-free}, we have:
\begin{propn}   \label{propn:euler-classifying}
Let $\scat{A}$ be a finite skeletal category containing no endomorphisms
except identities.  Then $\chi(B\scat{A})$ is defined and equal to
$\chi(\scat{A})$.  \done
\end{propn}

For the final compatibility result, consider the following schematic
diagrams:
\[
\begin{diagram}[size=1.5em]
\{ \textrm{triangulated manifolds} \}   &               &       \\
\dTo                                    &\rdTo>\chi     &       \\
\{ \textrm{posets} \}                   &\rTo_\chi      &\integers\\
\end{diagram}
\diagspace\diagspace
\begin{diagram}[size=1.5em]
\{ \textrm{triangulated orbifolds} \}   &               &       \\
\dTo                                    &\rdTo>\chi     &       \\
\{ \textrm{categories} \}               &\rTo_\chi      &\rationals.\\
\end{diagram}
\]
On the left, we start with a compact manifold $M$ equipped with a finite
triangulation.  As shown in~\S 3.8 of~\cite{Sta}, the topological Euler
characteristic of $M$ is equal to the Euler characteristic of the poset of
simplices in the triangulation, ordered by inclusion.  We generalize this
result from manifolds to orbifolds, which entails replacing posets by
categories and $\integers$ by $\rationals$.  

Let $M$ be a compact orbifold equipped with a finite triangulation.
(See~\cite{MP} for definitions.)  The simplices in the triangulation form a
poset $P$, and if $p \in P$ is a $d$-dimensional simplex then $\downset{p} =
\{ q \in P \such q \leq p \}$ is isomorphic to the poset $\scat{P}_{d + 1}$ of
nonempty subsets of $\{ 1, \ldots, d + 1 \}$, with $p \in \downset{p}$
corresponding to $\{ 1, \ldots, d + 1 \} \in \scat{P}_{d + 1}$.  Every $p \in
P$ has a stabilizer group $G(p)$, and
\[
\chi(M)
=
\sum_{p \in P} (-1)^{\dim p} / |G(p)|.
\]
On the other hand, the groups $G(p)$ fit together to form a \demph{complex of
finite groups} on $P^\op$, that is, a weak functor $G: P^\op \go \Cat$ taking
values in finite groups (regarded as one-object categories) and injective
homomorphisms; see~\S 3 of~\cite{Moe}.  This gives a finite category
$\elt{G}$.  For example, when $M$ is a manifold, each group $G(p)$ is trivial
and $\elt{G} \iso P$.  

The following result is joint with Ieke Moerdijk.

\begin{propn}
Let $M$ be a compact orbifold equipped with a finite triangulation.  Let $G$ be
the resulting complex of groups.  Then $\chi(\elt{G})$ is defined and equal
to $\chi(M)$.
\end{propn}

\begin{proof}
Every arrow in $\elt{G}$ is monic, so by Theorem~\ref{thm:mobius-idem-id},
$\elt{G}$ has Euler characteristic.  Moreover, $P$ is a finite poset, so has a
unique coweighting $k_\bl$, and $\chi(\elt{G}) = \sum_p k_p / |G(p)|$ by the
dual of Proposition~\ref{propn:euler-fib}.

The coweight of $p$ in $P$ is equal to the coweight of $p$ in $\downset{p}
\iso \scat{P}_{d + 1}$, where $d = \dim p$.  The unique coweighting $k_\bl$ on
$\scat{P}_{d + 1}$ is given by $k_J = (-1)^{|J| - 1}$, so $k_p = (-1)^{(d + 1)
- 1} = (-1)^{\dim p}$.  The result follows.  \done
\end{proof}

There is an accompanying theory of Lefschetz number.  Let $\scat{A}$ be a
finite category and $F: \scat{A} \go \scat{A}$ an endofunctor of
$\scat{A}$.  The category $\Fix{F}$ has as objects the (strict) fixed
points of $F$, that is, the objects $a \in \scat{A}$ such that $Fa = a$;
a map $a \go b$ in $\Fix{F}$ is a map $f: a \go b$ in $\scat{A}$ such that
$Ff = f$.  

\begin{defn}
Let $F$ be an endofunctor of a finite category.  Its \demph{Lefschetz
number} $\Lef{F}$ is $\chi(\Fix{F})$, when this exists.
\end{defn}

The Lefschetz number is, then, the sum of the (co)weights of the fixed points.
This is analogous to the standard Lefschetz fixed point formula, (co)weight
playing the role of index.  The following results further justify the
definition.

\begin{propn}
Let $\scat{A}$ be a finite category.
\begin{enumerate}
\item   \label{item:lef-id}
  $\Lef{1_{\scat{A}}} = \chi(\scat{A})$, one side being defined if and
  only if the other is.
\item   \label{item:lef-cycle}
  If $\scat{B}$ is another finite category and $\scat{A} \oppair{F}{G}
  \scat{B}$ are functors then $\Lef{GF} = \Lef{FG}$, one side being defined
  if and only if the other is.
\item   \label{item:lef-classifying}
  Let $F: \scat{A} \go \scat{A}$ and write $BF: B\scat{A} \go B\scat{A}$
  for the induced map on the classifying space of $\scat{A}$.  If
  $\scat{A}$ is skeletal and contains no endomorphisms except identities
  then $\Lef{F} = \Lef{BF}$, with both sides defined.
\end{enumerate}
\end{propn}
In the special case that $\scat{A}$ is a poset,
part~(\ref{item:lef-classifying}) is Theorem~1.1 of~\cite{BB}.

\begin{proof}
For~(\ref{item:lef-id}) and~(\ref{item:lef-cycle}), just note that
$\Fix{1_{\scat{A}}} \iso \scat{A}$ and $\Fix{GF} \iso \Fix{FG}$.
For~(\ref{item:lef-classifying}), recall from the proof of
Proposition~\ref{propn:euler-classifying} that $N\scat{A}$ has only
finitely many nondegenerate simplices; then
\begin{eqnarray*}
\Lef{BF}
        &=      &
\sum_{n \geq 0} 
(-1)^n
|\{ \textrm{nondegenerate $n$-simplices in $N\scat{A}$ fixed by $NF$} \}|
\\
        &=      &
\sum_{n \geq 0} 
(-1)^n
|\{ \textrm{nondegenerate $n$-paths in $\Fix{F}$} \}|
\\
        &=      &
\Lef{F},
\end{eqnarray*}
using Corollary~\ref{cor:mobius-circuit-free} in the last step.  
\done
\end{proof}

An \demph{algebra} for an endofunctor $F$ of $\scat{A}$ is an object
$a \in \scat{A}$ equipped with a map $h: Fa \go a$.  With the evident
structure-preserving morphisms, algebras for $F$ form a category $\Alg{F}$.
There is a dual notion of \demph{coalgebra} (where now $h: a \go Fa$), giving
a category $\Coalg{F}$.

\begin{propn}
Let $\scat{A}$ be a finite skeletal category containing no endomorphisms
except identities.  Then $\chi(\Alg{F}) = \Lef{F} = \chi(\Coalg{F})$, with all
three terms defined. 
\end{propn}

\begin{proof}
First observe that $\scat{A}$ is circuit-free.  Now, the inclusion $\Fix{F}
\go \Alg{F}$ has a right adjoint $R$: given an algebra $(a, h)$,
circuit-freeness implies that $F^n a$ is a fixed point for all sufficiently
large $n$, and $R(a, h) = F^n a$.  The Euler characteristics of $\Alg{F}$ and
$\Fix{F}$ exist, by Corollary~\ref{cor:mobius-circuit-free}, and are equal, by
Proposition~\ref{propn:euler-adjn-eqv}(\ref{item:euler-adjn}).  The statement
on coalgebras follows by duality.  \done
\end{proof}

For example, if $f$ is an endomorphism of a finite poset $A$ then the
sub-posets  
\[
\{ a \in A \such f(a) \leq a \},
\qquad
\{ a \in A \such f(a) = a \},
\qquad
\{ a \in A \such f(a) \geq a \}
\]
all have the same Euler characteristic. 

The theory of Euler characteristic presented here can be extended in at least
two directions.

First, we can relax the finiteness assumption.  For instance, the category
of finite sets and bijections should have Euler characteristic $\sum_{n =
0}^\infty 1/|S_n| = e$, as observed in~\cite{BD}.  See the remarks after
Corollary~\ref{cor:subcats}.

Second, note that the Euler characteristic of categories is defined in
terms of the cardinality of finite sets; it is then clear that the theory can
be developed for $\cat{V}$-enriched categories when there is a suitable notion
of cardinality or Euler characteristic of objects of $\cat{V}$.  For example,
$\cat{V}$ might be the category of finite-dimensional vector spaces, with
dimension playing the role of cardinality, and this leads to an Euler
characteristic for finite linear categories.  For another, recall that a
\demph{0-category} is a set and an $\demph{$n$-category}$ is a category
enriched in $(n - 1)$-categories; iterating, we obtain an Euler
characteristic for finite $n$-categories.  In particular, if $\mathbf{S}^n$ is
the $n$-category consisting of two parallel $n$-cells then
$\chi(\mathbf{S}^n) = 1 + (-1)^n$.

\section{The cardinality of a colimit}
\label{sec:colimit}

The main theorem of this section generalizes the formulas
\[
|X \cup Y| = |X| + |Y| - |X \cap Y|,
\diagspace\diagspace
|S/G| = |S|/|G|
\]
where $X$ and $Y$ are subsets of some larger set and $S$ is a set acted on
freely by a group $G$.  

Take a finite functor $X: \scat{A} \go \Set$.  The colimit $\colimit{\scat{A}}
X$ can be viewed as the gluing-together of the sets $Xa$.  Its cardinality
depends on the way in which these sets are glued together, which in turn is
determined by the action of $X$ on morphisms, so in general there is no
formula for $|\colimit{\scat{A}} X|$ purely in terms of the cardinalities
$|Xa|$ ($a \in \scat{A}$).

Suppose, however, that we are in the extreme case that there are no `unforced'
equations of the type $(Xf)(x) = (Xf')(x')$.  For pushouts, this means that
the two functions along which we are pushing out are injective; when
$\scat{A}$ is a group $G$, so that $X$ is a $G$-action, it means that the
action is free.  In this extreme case, $|\colimit{\scat{A}} X|$ can be
calculated as a weighted sum of the cardinalities $|Xa|$.

We now make this precise.  Recall that a $\Set$-valued functor is said to be
familially representable if it is a sum of representables.

\begin{propn}   \label{propn:colim-fr}
Let $\scat{A}$ be a finite category and $k^\bl$ a weighting on $\scat{A}$.
If $X: \scat{A} \go \Set$ is finite and familially representable then
$|\colimit{\scat{A}} X| = \sum_a k^a |Xa|$.
\end{propn}

\begin{proof}
The result holds if $X$ is representable, since then $|\colimit{\scat{A}} X| =
1$.  On the other hand, the class of functors $X$ for which the conclusion
holds is clearly closed under finite sums.  \done
\end{proof}

To make use of this, we need a way of recognizing familially representable
functors.  Carboni and Johnstone~\cite{CJ,CJcorr} show that when $\cat{A}$
satisfies certain hypotheses, including having all limits, a functor $\cat{A}
\go \Set$ is familially representable if and only if it preserves connected
limits.  This is of little help, because our categories $\scat{A}$ are finite,
and a finite category does not have even all finite limits unless it is a
lattice.

However, a standard philosophy applies: when $\scat{A}$ fails to have all
limits of a certain type, it is rarely useful to consider the functors
$\scat{A} \go \Set$ preserving limits of that type; the correct substitute is
the class of functors that are suitably `flat'.  The notion of flatness
appropriate here will be called nondegeneracy.  (This is unrelated to the
usage of `nondegenerate' in \S\ref{sec:mobius}.)

\begin{defn}    \label{defn:nondegen}
Let $\scat{A}$ be a small category.  A functor $X: \scat{A} \go \Set$ is
\demph{nondegenerate} if $\elt{X}$ has the following diagram-completion
properties:
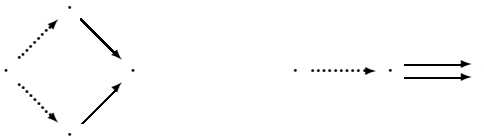
\[
\begin{diagram}[size=1.2em]
        &       &\cdot  &       &       \\
        &\ruGet &       &\rdTo  &       \\
\cdot   &       &       &       &\cdot  \\
        &\rdGet &       &\ruTo  &       \\
        &       &\cdot  &       &       \\
\end{diagram}
\diagspace \diagspace \diagspace
\begin{diagram}[size=1.8em]
\cdot   &\rGet  &\cdot  &\pile{\rTo\\ \rTo}     &\cdot  \\
\end{diagram}
\]
Explicitly, this means that (i) given arrows $a \goby{f} b \ogby{f'} a'$ in
$\scat{A}$ and $x \in Xa$, $x' \in Xa'$ satisfying $(Xf)(x) = (Xf')(x')$, there
exist arrows $a \ogby{g} c \goby{g'} a'$ and $z \in Xc$ satisfying $fg =
f'g'$, $(Xg)(z) = x$, and $(Xg')(z) = x'$, and (ii) given arrows $a
\parpair{f}{f'} b$ in $\scat{A}$ and $x \in Xa$ satisfying $(Xf)(x) =
(Xf')(x)$, there exist $c \goby{g} a$ and $z \in Xc$ satisfying $fg = f'g$
and $(Xg)(z) = x$.  
\end{defn}

This is the most concrete form of the definition.  For further explanation
and a proof that nondegeneracy is equivalent to familial representability, see
the Appendix; for references, see~\cite{SS1}. 

From Lemma~\ref{lemma:fr-nondegen} we deduce a more applicable form of
Proposition~\ref{propn:colim-fr}:

\begin{thm} \label{thm:colim-nd}
Let $\scat{A}$ be a finite Cauchy-complete category and $k^\bl$ a weighting
on $\scat{A}$.  If $X: \scat{A} \go \Set$ is finite and nondegenerate then
$|\colimit{\scat{A}} X| = \sum_a k^a |Xa|$.
\done
\end{thm}

Recalling that $\colimit{\scat{A}} X$ is the set of connected-components of
$\elt{X}$, this may be rephrased as $\pi_0(\elt{X}) = \sum k^a |Xa|$.  On the
other hand, Proposition~\ref{propn:euler-fib} implies that $\chi(\elt{X}) =
\sum k^a |Xa|$.  Indeed, under the hypotheses of the Theorem, $X$ is
familially representable, so each connected-component of $\elt{X}$ has an
initial object, so $\chi(\elt{X}) = \pi_0(\elt{X})$.

\begin{examples}    \label{egs:colim-formula}
\begin{enumerate}
\item \label{eg:colim-pushout} Let $\scat{L}$ be the category
  of~\ref{egs:wtgs}(\ref{eg:wtgs-pushout}).  A functor $X: \scat{L} \go
  \Set$ is nondegenerate if and only if both functions $Xa \go Xb_i$ are
  injective.  In that case, Theorem~\ref{thm:colim-nd} says that
\[
| Xb_1 +_{Xa} Xb_2 |
=
|Xb_1| + |Xb_2| - |Xa|
\]
where the set on the left-hand side is a pushout.

\item   \label{eg:colim-coeq}
Let $\scat{B}$ be the category $\left( a \parpair{f}{g} b \right)$.  A functor
$X: \scat{B} \go \Set$ is nondegenerate if and only if the two functions $Xf$,
$Xg$ are injective and have disjoint images.  The unique weighting $k^\bl$ on
$\scat{B}$ is $(k^a, k^b) = (-1, 1)$, and
\[
| (Xb)/\sim |
=
|Xb| - |Xa|
\]
where $\sim$ is the equivalence relation generated by $(Xf)(x) \sim
(Xg)(x)$ for all $x \in Xa$.

\item Let $G$ be a group.  A functor $X: G \go \Set$ is a set $S$ equipped
with a left $G$-action; the functor is nondegenerate if and only if the action
is free.  Theorem~\ref{thm:colim-nd} then says that the number of orbits is
$|S|/|G|$.

\item The Theorem can be viewed as a generalized inclusion-exclusion
principle.  (Compare~\cite{Rota}.)  Let $n \geq 0$ and let $\scat{P}_n$ be the
poset of nonempty subsets of $\{1, \ldots, n\}$, ordered by inclusion. (So
$\scat{P}_2^\op$ is the category $\scat{L}$ of~(\ref{eg:colim-pushout}).)  Its
unique coweighting $k_\bl$ is defined by $k_J = (-1)^{|J| - 1}$.  Given
subsets $S_1, \ldots, S_n$ of some set, there is a nondegenerate functor $X:
\scat{P}_n^\op \go \Set$ defined on objects by $X(J) = \bigcap_{j \in J} S_j$
and on maps by inclusion.  Theorem~\ref{thm:colim-nd} gives the
inclusion-exclusion formula,
\[
|S_1 \cup \cdots \cup S_n|
=
\sum_{\emptyset \neq J \sub \{ 1, \ldots, n \} }
(-1)^{|J| - 1}
\left|
\bigcap_{j \in J} S_j
\right|.
\]
\end{enumerate}
\end{examples}

\begin{cor}
Let $\scat{A}$ be a finite Cauchy-complete category admitting a weighting.
Let $X, Y: \scat{A} \go \Set$ be finite nondegenerate functors satisfying
$|Xa| = |Ya|$ for all $a \in \scat{A}$.  Then $|\colimit{\scat{A}} X| =
|\colimit{\scat{A}} Y|$.  \done
\end{cor}

The condition that $\scat{A}$ admits a weighting cannot be dropped:
consider the category $\scat{A}$ of
Example~\ref{egs:wtgs}(\ref{eg:wtgs-none}) and the functors $X =
\scat{A}(a_1, \dashbk) + \scat{A}(a_4, \dashbk)$, $Y = \scat{A}(a_2,
\dashbk)$.

If $\scat{A}$ not only has a weighting but admits M\"obius inversion then
a stronger statement can be made: Proposition~\ref{propn:fr-coeffs}.

\section{Relations with Rota's theory}
\label{sec:other}

In 1964, Gian-Carlo Rota published his seminal paper~\cite{Rota} on M\"obius
inversion in posets.  The name is motivated as follows: in the poset of
positive integers ordered by divisibility, $\mu(a, b) = \mu(b/a)$ whenever $a$
divides $b$, where the $\mu$ on the right-hand side is the classical M\"obius
function.  He was not the first to define M\"obius inversion in
posets---Weisner, Hall, and Ward preceded him---but Rota's contribution was
the decisive one; in particular, he realized the power of the method in
enumerative combinatorics.  The history of M\"obius inversion is well
described in~\cite{Rota}, \cite{Gre} and~\cite{Sta}.

In this section we see that some of the principal results in Rota's theory are
the order-theoretic shadows of more general categorical facts.  We also
consider a different generalization of M\"obius--Rota inversion.

Given a poset $A$, Rota considered its \demph{incidence algebra} $I(A)$, which
is the subring of $\inc{A}$ consisting of the integer-valued $\theta \in
\inc{A}$ such that $\theta(a, b) = 0$ whenever $a \not\leq b$.  By
Example~\ref{egs:mobius}(\ref{eg:mobius-posets}) or
Corollary~\ref{cor:mobius-circuit-free}, $\mu \in I(A)$.

In posets, then, $\zeta(a, b) = 0 \implies \mu(a, b) = 0$.  More generally:
\begin{thm} \label{thm:mobius-zero}
If $\scat{A}$ is a finite category with M\"obius inversion then, for $a, b
\in \scat{A}$,
\[
\zeta(a, b) = 0 \implies \mu(a, b) = 0.
\]
\end{thm}

The proof uses a combinatorial lemma.

\begin{lemma}
Let $n \geq 2$ and $\sigma \in S_{n - 1}$.  Then there exist $k \geq 1$ and
$p_0, \ldots, p_k$ such that
\[
p_0 = 1,
\quad
p_1, \ldots, p_{k - 1} \in \{ 1, \ldots, n - 1 \},
\quad
p_k = n,
\]
and $p_r = \sigma(p_{r - 1}) + 1$ for each $r \in \{1, \ldots, k \}$.
\end{lemma}

\begin{proof}
Suppose not; then there is an infinite sequence $(p_r)_{r \geq 0}$ of elements
of $\{ 1, \ldots, n - 1 \}$ satisfying $p_0 = 1$ and $p_r = \sigma(p_{r - 1})
+ 1$ for all $r \geq 1$.  Let $\epsln$ be the endomorphism of the finite set
$\{ p_r \such r \geq 0 \}$ defined by $\epsln(p) = \sigma(p) + 1$.  Then
$\epsln$ is injective but not surjective (since $1 \not\in \im\,\epsln$),
contradicting finiteness.  \done
\end{proof}

\begin{prooflike}{Proof of Theorem~\ref{thm:mobius-zero}}
Write the objects of $\scat{A}$ as $a_1, \ldots, a_n$.  There is an $n\times
n$ matrix $Z$ defined by $Z_{ij} = \zeta(a_i, a_j)$, and $Z$ is invertible
over $\rationals$ with $(Z^{-1})_{ij} = \mu(a_i, a_j)$.  Suppose that $i, j
\in \{1, \ldots, n\}$ and $Z_{ij} = 0$. Certainly $i \neq j$, so $n \geq 2$
and we may assume that $(i, j) = (1, n)$.  By Cramer's formula for the inverse
of a matrix, our task is to prove that the $(n, 1)$-minor $Z'$ of $Z$ has
determinant zero.

The $(p, q)$-entry of $Z'$ is $Z_{p, q + 1}$, so
\[
\det Z'
=
\sum_{\sigma \in S_{n - 1}}
\pm Z_{1, \sigma(1) + 1} \cdots Z_{n - 1, \sigma(n - 1) + 1}.
\]
It suffices to prove that each summand is $0$.  Indeed, let $\sigma \in
S_{n - 1}$.  Take $p_0, \ldots, p_k$ as in the Lemma.  By hypothesis, there
is no map $a_1 \go a_n$ in $\scat{A}$.  Categories have composition, so
there is no diagram
\[
a_1 = a_{p_0} \go a_{p_1} \go \ \cdots \ \go a_{p_k} = a_n
\]
in $\scat{A}$.  Hence $\zeta(a_{p_{r - 1}}, a_{p_r}) = 0$ for some $r \in
\{ 1, \ldots, k \}$, giving $Z_{p_{r - 1}, \sigma(p_{r - 1}) + 1} = 0$,
as required.  \done
\end{prooflike}

Given objects $a, c$ of a category $\scat{A}$, let $\scat{A}_{a, c}$ be the
full subcategory consisting of those $b \in \scat{A}$ for which there exist
arrows $a \go b \go c$.  Theorem~\ref{thm:mobius-zero} easily implies:

\begin{cor}     \label{cor:subcats}
Let $\scat{A}$ be a finite category.  Then $\scat{A}$ has M\"obius inversion
if and only if $\scat{A}_{a, c}$ has M\"obius inversion for all $a, c \in
\scat{A}$, and in that case the M\"obius function of $\scat{A}_{a, c}$ is the
restriction of that of $\scat{A}$.
\done
\end{cor}

These results suggest a way of relaxing the finiteness assumption on our
categories.  It extends to categories the local finiteness condition on posets
used in the Rota theory.  Let $\scat{A}$ be a category for which each
subcategory $\scat{A}_{a, c}$ is finite.  Then each hom-set $\scat{A}(a, b)$
has finite cardinality, $\zeta(a, b)$, and there is a $\rationals$-algebra
\[
\incinf{\scat{A}}
=
\{  
\theta: \ob\scat{A} \times \ob\scat{A} \go \rationals
\such
\textrm{for $a, b \in \scat{A}$, $\zeta(a, b) = 0 \implies \theta(a, b) =
0$} 
\}
\]
with operations defined as for $\inc{\scat{A}}$.  Evidently $\zeta \in
\incinf{\scat{A}}$, and $\scat{A}$ may be said to have M\"obius inversion if
$\zeta$ has an inverse $\mu$ in $\incinf{\scat{A}}$.  By
Theorem~\ref{thm:mobius-zero}, this extends the definition for finite
categories.  For example, the skeletal category $\Dinj{}$ of finite totally
ordered sets and order-preserving injections has M\"obius inversion; compare
Example~\ref{egs:mobius}(\ref{eg:mobius-Dinj}).

The main theorem in Rota's paper~\cite{Rota} relates the M\"obius functions of
two posets linked by a Galois connection.  Viewing a poset as a special
category, a Galois connection is nothing but a (contravariant) adjunction,
which suggests the following generalization of Rota's theorem.
\begin{propn}
Let $\scat{A}$ and $\scat{B}$ be finite categories with M\"obius inversion.
Let $\scat{A} \oppair{F}{G} \scat{B}$ be an adjunction, $F \ladj G$.  Then
for all $a \in \scat{A}$, $b \in \scat{B}$,
\[
\sum_{a': Fa' = b} \mu(a, a')
=
\sum_{b': Gb' = a} \mu(b', b).
\]
\end{propn}

\begin{proof}
Write $\zeta(a, b) = \zeta(Fa, b) = \zeta(a, Gb)$.  Then for all $a \in
\scat{A}$, $b \in \scat{B}$,
\[
\sum_{a': Fa' = b} \mu(a, a')   
=       
\sum_{a' \in \scat{A}} \mu(a, a') \delta(Fa', b)
=       
\sum_{a' \in \scat{A},\, b' \in \scat{B}}
\mu(a, a') \zeta(a', b') \mu(b', b). 
\]
The result follows by symmetry.
\done
\end{proof}

For example, when $l$ is an element of a finite lattice $L$, the inclusion of
the sub-poset $\{ x \in L \such x \leq l \}$ into $L$ has right adjoint
$(\dashbk \wedge l)$, giving Weisner's Theorem (p.351 of~\cite{Rota}).

The Euler characteristic of posets has been studied extensively;
see~\cite{Sta} for references.  Given a finite poset $A$, the classifying
space $BA$ always has Euler characteristic, which by
Proposition~\ref{propn:euler-classifying} is equal to the Euler characteristic
of the category $A$.  On the other hand, we may form a new poset $\twid{A}$ by
adjoining to $A$ a least element $0$ and a greatest element $1$, and then
$\chi(A) = \mu_{\twid{A}}(0, 1) + 1$; see~\cite{Rota} or~\S 3.8 of~\cite{Sta}.
This result can be extended from posets to categories:
\begin{propn}
Let $\scat{A}$ be a finite category.  Write $\twid{\scat{A}}$ for the
category obtained from $\scat{A}$ by freely adjoining an initial object $0$
and a terminal object $1$.  If $\scat{A}$ has M\"obius inversion then
$\twid{\scat{A}}$ does too, and $\mu_{\twid{\scat{A}}}(0, 1) =
\chi(\scat{A}) - 1$.
\end{propn}

\begin{proof}
Suppose that $\scat{A}$ has M\"obius inversion.  Let $\scat{A}_0$ be the
category obtained from $\scat{A}$ by freely adjoining an initial object
$0$.  Extend $\mu \in \inc{\scat{A}}$ to a function $\mu \in
\inc{\scat{A}_0}$ by defining
\[
\mu(0, b) = - \sum_{a \in \scat{A}} \mu(a, b),
\qquad
\mu(a, 0) = 0,
\qquad
\mu(0, 0) = 1
\]
($b, a \in \scat{A}$).  It is easily checked that this is the
M\"obius function of $\scat{A}_0$.

Dually, if $\scat{B}$ is a finite category with M\"obius inversion then the
category $\scat{B}_1$ obtained from $\scat{B}$ by freely adjoining a
terminal object $1$ also has M\"obius inversion, with $\mu(c, 1) = -
\sum_{b \in \scat{B}} \mu(c, b)$ for all $c \in \scat{B}$.  Take
$\scat{B} = \scat{A}_0$: then $\scat{A}_{01} = \twid{\scat{A}}$ has
M\"obius inversion, and
\[
\mu(0, 1) 
=
-\sum_{b \in \scat{A}_0} \mu(0, b)
=
-\sum_{b \in \scat{A}} \mu(0, b) - \mu(0, 0)
=
\sum_{a, b \in \scat{A}} \mu(a, b) - 1
=
\chi(\scat{A}) - 1.
\]
\done
\end{proof}

\paragraph*{Remark} Recall~\cite{CKW} that given categories $\scat{B},
\scat{A}$ and a functor $M: \scat{B}^\op \times \scat{A} \go \Set$, the
\demph{collage} of $M$ is the category $\scat{C}$ formed by taking the
disjoint union of $\scat{B}$ and $\scat{A}$ and adjoining one arrow $b \go
a$ for each $b \in \scat{B}$, $a \in \scat{A}$ and $m \in M(b, a)$, with
composition defined using $M$.  Assuming finiteness, if $\scat{B}$ and
$\scat{A}$ have M\"obius inversion then so does $\scat{C}$:
\[
\renewcommand{\arraystretch}{1.5}
\begin{array}{l}
\displaystyle
\mu_{\scat{C}}(b, b') = \mu_{\scat{B}}(b, b'),
\qquad
\mu_{\scat{C}}(a, a') = \mu_{\scat{A}}(a, a'),
\qquad
\mu_{\scat{C}}(a, b) = \emptyset,\\
\displaystyle
\mu_{\scat{C}}(b, a)
=
- \sum_{b', a'} 
\mu_{\scat{B}}(b, b')\, |M(b', a')|\, \mu_{\scat{A}}(a', a)
\end{array}
\renewcommand{\arraystretch}{1}
\]
($b, b' \in \scat{B}$, $a, a' \in \scat{A}$).  In the proof above, the
calculation of the M\"obius function of $\scat{A}_0$ is the special case where
$\scat{B}$ is the terminal category and $M$ has constant value $1$.  The
ordinal sum of posets is another special case.  Moreover, one easily deduces a
formula for the Euler characteristic of a collage, which in the special case
of posets is essentially Theorem~3.1 of Walker~\cite{Walk}.

\paragraph*{}
Let us now look at the different generalization of Rota's M\"obius inversion
proposed, independently, by Content, Lemay and Leroux~\cite{CLL} and
by Haigh~\cite{Hai}.  (See also~\cite{Ler} and~\S 4 of~\cite{Law}.  Haigh
briefly considered the same generalization as here, too; see~3.5
of~\cite{Hai}.)  Given a sufficiently finite category $\scat{A}$, they take
the algebra $I(\scat{A})$ of functions from $\{\textrm{arrows of }\scat{A}\}$
to $\rationals$ (or more generally, to some base commutative ring), with a
convolution product:
\[
(\theta\phi)(f) 
=
\sum_{hg = f} \theta(g) \phi(h).
\]
Taking $\zeta \in I(\scat{A})$ to have constant value $1$, they call the
\demph{M\"obius function} of $\scat{A}$ the inverse $\mu = \zeta^{-1}$ in
$I(\scat{A})$, if it exists.  When $\scat{A}$ is a poset, this agrees with
Rota; when $\scat{A}$ is a monoid, it agrees with Cartier and Foata~\cite{CF}.

They seek to solve a harder problem than we do: if a finite category
$\scat{A}$ has M\"obius inversion in their sense then it does in ours (with
$\mu(a, b) = \sum_{f \in \scat{A}(a, b)} \mu(f)$), but not conversely.  For
instance, a non-trivial finite group never has M\"obius inversion in their
sense, but always does in ours.

\section{Appendix: category theory}
\label{sec:app}

Here follows a skeletal account of some standard notions: category of
elements, flat functors, and Cauchy-completeness.  Details can be found in
texts such as~\cite{Bor}.  Throughout, $\scat{A}$ denotes a small
category.  

Let $X: \scat{A} \go \Set$.  The \demph{category of elements} $\elt{X}$ of $X$
has as objects all pairs $(a, x)$ where $a \in \scat{A}$ and $x \in Xa$, and
as maps $(a, x) \go (a', x')$ all maps $f: a \go a'$ in $\scat{A}$ such that
$(Xf)(x) = x'$.

Similarly, if $X: \scat{A} \go \Cat$ (where $\Cat$ is the category of small
categories and functors) then $X$ has a \demph{category of elements}
$\elt{X}$; its objects are pairs $(a, x)$ where $a \in \scat{A}$ and $x \in
Xa$, and its maps $(a, x) \go (a', x')$ are pairs $(f, \xi)$ where $f: a \go
a'$ in $\scat{A}$ and $\xi: (Xf)(x) \go x'$ in $Xa'$.  This definition can be
made even when $X$ is a \demph{weak functor} or \demph{pseudofunctor}, that
is, only preserves composition and identities up to coherent isomorphism.  The
weak functors $\scat{A} \go \Cat$ correspond to the fibrations over
$\scat{A}^\op$; see~\cite{Bor}.

The definition for $\Cat$-valued functors extends that for $\Set$-valued
functors if a set is viewed as a discrete category (one with no maps other
than identities).

Any two functors $Y: \scat{A}^\op \go \Set$ and $X: \scat{A} \go \Set$ have a
tensor product $Y \otimes X$, a set, defined by
\[
Y \otimes X 
=
\left( \coprod_{a \in \scat{A}} Ya \times Xa \right) / \sim
\]
where $\sim$ is the equivalence relation generated by $(y, (Xf)(x)) \sim
((Yf)(y), x)$ whenever $f: a \go b$, $x \in Xa$ and $y \in Yb$.  (It may be
helpful to think of $X$ and $Y$ as left and right $\scat{A}$-modules.)  A
functor $X: \scat{A} \go \Set$ is \demph{flat} if
\[
\dashbk \otimes X: \ftrcat{\scat{A}^\op}{\Set} \go \Set
\]
preserves finite limits.  An equivalent condition is that $\elt{X}$ is
\demph{cofiltered}, that is, every finite diagram in $\elt{X}$ admits at least
one cone.  

\begin{propn}   \label{propn:nondegen}
The following conditions on a functor $X: \scat{A} \go \Set$ are equivalent:
\begin{enumerate}
\item $X$ is nondegenerate (in the sense of~\ref{defn:nondegen})
\item every connected-component of $\elt{X}$ is cofiltered
\item $X$ is a sum of flat functors.
\item $\dashbk \otimes X: \ftrcat{\scat{A}^\op}{\Set} \go \Set$ preserves
finite connected limits
\end{enumerate}
\end{propn}

\begin{proof}
See~\cite{SS1} or~\cite{ABLR}.
\done
\end{proof}

An idempotent $e: a \go a$ in $\scat{A}$ \demph{splits} if there exist $a
\oppair{s}{i} b$ such that $si = 1$ and $is = e$.  The category $\scat{A}$ is
\demph{Cauchy-complete} if every idempotent in $\scat{A}$ splits.
All of the examples of categories in this paper are Cauchy-complete, except
that a finite monoid is Cauchy-complete if and only if it is a group.

\begin{lemma}   \label{lemma:fr-nondegen}
Let $\scat{A}$ be a Cauchy-complete category and $X: \scat{A} \go \Set$ a
finite functor.  Then $X$ is familially representable if and only if $X$ is
nondegenerate.
\end{lemma}

\begin{proof}
By Proposition~\ref{propn:nondegen}, it is enough to prove that a finite
functor $X$ is representable if and only if it is flat.  `Only if' is
immediate.

For `if', suppose that $X$ is flat.  Then $\elt{X}$ is cofiltered, and finite
by hypothesis, so the identity functor $1_{\elt{X}}$ admits a cone.  Also,
$\elt{X}$ is Cauchy-complete since $\scat{A}$ is.  Now, if $\scat{C}$ is a
Cauchy-complete category and $(j \goby{p_c} c)_{c \in \scat{C}}$ is a cone on
$1_{\scat{C}}$ then $p_j$ is idempotent, and the object through which it
splits is initial.  Hence $\elt{X}$ has an initial object; equivalently, $X$
is representable.  \done
\end{proof}

\end{document}